\documentclass{article}
\usepackage{amsmath,amssymb,amsfonts, amsthm}
\usepackage[dvipsnames]{xcolor}
\usepackage{graphicx} 
\usepackage{tikz}
\usepackage{tkz-base}
\usepackage{tkz-euclide}

\usepackage[colorlinks=true, pdfstartview=FitV, linkcolor=blue, citecolor=blue, urlcolor=blue]{hyperref}
\usepackage[colorinlistoftodos]{todonotes}
\usepackage{dsfont}
\usepackage[russian,english]{babel}

\definecolor{darkred}{rgb}{0.7,0,0} 
\newcommand{\dmu}{\; \mathrm{d\mu}}
\newcommand{\dx}{\,\mathrm{dx}}
\newcommand{\dy}{\,\mathrm{dy}}
\newcommand{\supp}{\mathrm{supp}}

\newcommand{\CC}{\mathbb{C}}
\newcommand{\GL}{\mathrm{GL}}
\newcommand{\tr}{\mathrm{tr}}
\newcommand{\JUE}{\mathrm{JUE}}
\newcommand{\Howe}{\mathrm{H}}

\author{Anton Andreevich Nazarov$^{1,2,3}$, Matvey Stanislavovich Sushkov$^{4,5}$\\
{\small$^{1}$Department of Physics, St. Petersburg State University,} \\
{\small  Ulyanovskaya 1, 198504 St.~Petersburg, Russia}\\
{\small$^{2}$ Beijing Institute of Mathematical Sciences and Applications (BIMSA),}\\
{\small Bejing 101408, People’s Republic of China}\\
{\small$^{3}$ email:antonnaz@gmail.com}\\
{\small$^{4}$ 
  St. Petersburg Department
    of Steklov Mathematical Institute}\\
  {\small of Russian Academy of Sciences}\\
  {\small 191023, Fontanka 27, St. Petersburg, Russia}\\
{\small$^{5}$ email: matveysushkov34@gmail.com}
}

\title{Random Young diagrams and Jacobi Unitary Ensemble}

\date{October 2025}

\begin{document}

\maketitle
\begin{abstract}
  We consider random Young diagrams with respect to the measure induced by the decomposition of the $p$-th exterior power of $\mathbb{C}^{n}\otimes \mathbb{C}^{k}$ into irreducible representations of $\GL_{n}\times\GL_{k}$. We demonstrate that transition probabilities for these diagrams in the limit $n,k,p\to\infty$ with $p\sim nk$ converge to the large $N$ limiting law for the eigenvalues of random matrices in Jacobi Unitary Ensemble. We compute the characters of Young--Jucys--Murphy elements in $\bigwedge^{p}(\mathbb{C}^{n}\otimes\mathbb{C}^{k})$ and discuss their relation to surface counting. We formulate several conjectures on the connection  between the correlators in both random ensembles. 
\end{abstract}

\section{Introduction}

Transition probability for random Young diagrams attained prominent role in the works of S.V.~Kerov~\cite{kerov1993transition}. It was originally defined on the Young graph which is graded graph that has Young diagrams $\nu\vdash n$ of $n$ boxes on $n$-th level that are connected to the diagrams $\lambda$ of $n+1$ box if $\lambda$ is obtained from $\nu$ by addition of a single box, $\nu\nearrow\lambda$. Transition probability
\begin{equation}
  \label{eq:transition-probability}
  P(\nu\nearrow\lambda)=\frac{\dim\lambda}{|\lambda|\dim\nu}
\end{equation}
is then naturally related to the Plancherel measure on Young diagrams 
\begin{equation*}
  \label{eq:plancherel-measure}
  \mu^{\mathrm{Pl}}(\lambda)=\frac{(\dim \lambda)^{2}}{|\lambda|!}
\end{equation*}
by conditional probability formula
\begin{equation*}
  \label{eq:conditional-probability}
 \mu^{\mathrm{Pl}}(\lambda)=\sum_{\nu\nearrow\lambda}P(\nu\nearrow\lambda)\mu^{\mathrm{Pl}}(\nu).
\end{equation*}
Transition probability can be seen as a discrete measure defined on the ``corner boxes'' of the diagram $\lambda$ that could be removed to produce diagrams $\nu$. Similarly, cotransition probability can be introduced on the positions where a corner box can be added.  We will identify Young diagram with the sequence of particles and holes on the half-integer lattice $\mathbb{Z}+\frac{1}{2}$ or Maya diagram by rotating Young diagram in French notation $45^{\circ}$ and putting a particle in the midpoint of any unit interval on which the upper boundary of the rotated diagram decreases. Then transition probability becomes a discrete measure on half-integer lattice.

For example, Maya diagram for $\lambda=(5,5,2,1)$ looks as follows, transition probability is defined on shaded boxes or on corresponding particle positions that are marked by the larger circles:
\[
\begin{tikzpicture}[scale=.4,>=latex]
\draw[-,very thick] (0,0) -- ++(6.2,6.2);
\draw[-,very thick] (0,0) -- ++(-6.2,6.2);
\draw[<->,thick] (-6.5,0) -- (6.5,0);
\foreach \la/\i in {4/1,3/2,2/3,2/4,2/5}
  \draw[-] (\i,\i) -- ++(-\la,\la);
\foreach \la/\i in {5/1,5/2,2/3,1/4} {
  \draw[-] (-\i,\i) -- ++(\la,\la);
  \draw[color=blue,dashed] (-\i,\i) + (\la+.5,\la-.5) -- (\la-\i+.5,0);
  \fill[color=darkred] (-\i,\i) + (\la+.5,\la-.5) circle (0.1);
  \fill[color=darkred] (\la-\i+.5,0) circle (0.1);
}
\foreach \i in {5,6} {
  \draw[-] (-\i,\i) -- ++(.1,.1);
  \draw[color=blue,dashed] (-\i+.5,\i-.5) -- (-\i+.5,0);
  \fill[color=darkred] (-\i+.5,\i-.5) circle (0.1);
  \fill[color=darkred] (-\i+.5,0) circle (0.1);
}
\foreach \i in {-6,-5,...,6} {
  \draw[-] (\i,0.2) -- (\i,-0.2) node[anchor=north] {{\tiny$\i$}};
}
\foreach \i/\h in {-3/3,-1/3,3/5} {
  \path[fill=black,opacity=0.2] (\i,\h) -- (\i+1,\h+1) -- (\i,\h+2) -- (\i-1,\h+1);
  \fill[color=darkred] (\i+.5,0) circle (0.18);
  }
\end{tikzpicture}
\]

S.V.~Kerov had shown that in the limit $n\to\infty$ and scaling lattice step by $\frac{1}{\sqrt{n}}$, the transition measure for Plancherel-random Young diagrams converges to the probability measure with the density $\sigma(x)=\frac{\sqrt{4-x^{2}}}{2\pi}$ which is the large $N$ limiting density of eigenvalues of the Gaussian Unitary Ensemble (GUE) of random $N\times N$ hermitian matrices (Wigner's semicircle law) \cite{kerov1993transition}, see also \cite{kerov2003asymptotic}. In the same limit random Young diagrams converge in probability in supremum norm to the celebrated Vershik--Kerov--Logan--Schepp limit shape $\Omega(x)$ \cite{vershik1977asymptotics,logan1977variational} which is related to the limiting density of particles $\rho(x)=\arcsin\left(\frac{x}{2}\right)$ by $\Omega'(x)=1-2\rho(x)$.

Moreover, the bijection between the diagram profile $\Omega(x)$ and the transition probability measure $\sigma$ was termed Markov--Krein correspondence. It can be stated as a relation between two measures without relying on the limiting procedure as
\begin{equation*}
  \label{eq:markov-krein-profile}
  \int_a^b\frac{\mathrm{d}\sigma(x)}{u-x}=\frac{1}{u}\exp\left(\frac12 \int_a^b\frac{\mathrm{d}\left(\Omega(x)-|x|\right)}{x-u}\right).
\end{equation*}
where $a$ and $b$ denote the endpoints of the diagram, $(a,b) = \supp(\Omega(x)-|x|)$. This correspondence can also be expressed in terms of the particle density $\rho(x) = \frac{1 - \Omega'(x)}{2}$ with $\operatorname{supp}\rho \subset (-\infty, b]$ as
\begin{equation}
  \label{eq:markov-krein}
  \int_a^b\frac{\mathrm{d}\sigma(x)}{u-x}=\frac{1}{u-a}\exp\left( \int_a^b\frac{\rho(x)\mathrm{dx}}{u-x}\right).
\end{equation}

Later the same relation was demonstrated to hold between the limit shape obtained by P.~Biane for Schur--Weyl duality and Marchenko--Pastur limiting law for Laguerre Unitary Ensemble (LUE) or Wishart random matrices  \cite{biane2001approximate}. It is possible to extend this correspondence by constructing a piecewise-linear function $f$ by interlacing eigenvalues of random matrix $X$ from GUE or LUE with eigenvalues of $\hat X$ obtained from $X$ by removing last row and last column and then taking eigenvalues of $X$ to be positions of minima and eigenvalues of $\hat X$ to be positions of maxima and requiring $f'(x)=1$ between a minimum and the next maximum and $f'(x)=-1$ between a maximum and the next minimum. Such a function for GUE in the limit $N\to\infty$ converges to Vershik--Kerov--Logan--Schepp limit shape and to Biane's limit shape for LUE \cite{meliot2011kerov,bufetov2013kerov}.

The third classical random matrix ensemble which is known as Jacobi Unitary Ensemble (JUE) \cite{forrester2010log} or MANOVA matrix distribution \cite{dumitriua2002matrix} seems to be missing in the literature from this discussion of connection to random Young diagrams.

In the present paper we demonstrate that the limiting law for JUE is the limit of transition probability of random Young diagrams with respect to the measure induced by skew $(\GL_{n},\GL_{k})$ Howe duality~\cite{howe1989remarks}. We consider the space $\bigwedge^{p}(\mathbb{C}^{n}\otimes\mathbb{C}^{k})$, which has commuting actions of $\GL_{n}$ and $\GL_{k}$ and decomposes multiplicity-free into irreducible $\GL_{n}\times\GL_{k}$-representations. Taking the ratio of the dimension of the irreducible component corresponding to Young diagram $\lambda$ to the dimension of the whole space we get the probability measure on Young diagrams of $p$ boxes inside of $n\times k$ rectangle
\begin{equation}
  \label{eq:skew-howe-probability-measure}
  \mu^{\Howe}(\lambda)=\frac{\dim V_{\GL_{n}}(\lambda) \dim V_{\GL_{k}}(\lambda')}{\dim \bigwedge^{p}(\mathbb{C}^{n}\otimes\mathbb{C}^{k})}.
\end{equation}
In the limit $n,k,p\to\infty$ with $p\sim nk$ the convergence to the limit shape was established several times by different methods and different authors, see for example \cite{ismail1998strong,pittel2007limit,GTW01}. In Section \ref{sec:limiting-law-jue} we demonstrate that the limit shape of random Young diagrams with respect to $\mu^{\Howe}$ is connected to the limiting law for the eigenvalues of JUE by Markov--Krein correspondence.

We expect that similar correspondence takes place for all correlation functions as well because of a deeper structure. Okounkov in \cite{okounkov2000random} has demonstrated that correlators of row lengths of random Young diagrams with respect to Plancherel measure are traces of powers of Young--Jucys--Murphy elements. He also has shown that these traces count the number of certain Riemann surfaces. On the other hand GUE correlators also count surfaces. Okounkov proved that, asymptotically, these two sets of surfaces coincide. We hope to prove the connection of JUE correlators and correlators of row lengths of random diagrams with respect to $\mu^{\Howe}$ in a similar way. Therefore in  Section \ref{sec:rand-diagr-count} we compute the characters of the powers of Young--Jucys--Murphy  elements on $\bigwedge^{p}(\mathbb{C}^{n}\otimes\mathbb{C}^{k})$.  In the finite $n,k,p$ setting we observe the appearance of cotransition probabilities. We then explain that these characters count the number of classes of certain Riemann surfaces. We also recall existing results that connect JUE correlation functions to Hurwitz numbers that also count Riemann surfaces and state the problem of asymptotic correspondence.  

In Conclusion we state three conjectures. First is that fluctuations of random Young diagrams and JUE random matrices correspond to the same asymptotic set of Riemann surfaces. This conjecture is supported by the fact that fluctuations in both ensembles are described by the Airy kernel, and by the analogous statement for GUE and Plancherel-random diagrams proven by A.~Okounkov~\cite{okounkov2000random}. Second, the generating function for JUE correlators is tau function of an integrable hierarchy and admits a free-fermionic formulation. There is a free-fermionic construction for the correlation kernel for $\mu^{\Howe}$. We conjecture that these free-fermionic constructions are asymptotically connected. At last, we conjecture that piecewise-linear function obtained by interlacing JUE eigenvalues as described above converges to the limit shape for $\mu^{\Howe}$ and present numerical evidence of that. 

\section*{Acknowledgements}
The authors thank Alexander Tikhomirov for the references \cite{pastur2000law,ledoux2004differential}, Pavel Nikitin for useful discussions and Aleksandr Orlov for useful comments and references on tau-functions. 
The work of Matvey Sushkov was supported by the Ministry of Science and Higher Education of the Russian Federation (agreement 075-15-2025-344 dated 29/04/2025 for Saint Petersburg Leonhard Euler International Mathematical Institute at PDMI RAS).

\section{Limiting law for JUE eigenvalues and limit shape of Young diagrams}
\label{sec:limiting-law-jue}

The Jacobi Unitary Ensemble (JUE) can be defined as a probability measure on the set of hermitian matrices $H$ of size $N\times N$ with eigenvalues in the interval~$(0,1)$:
\begin{equation*}
  \label{eq:jacobi-ensemble}
  \mathrm{dm}^{\JUE}_{N}(H)=\frac{1}{Z^{\JUE}} (\det H)^{\alpha} (\det(I-H))^{\beta} \mathrm{d H},
\end{equation*}
where $\mathrm{d H}= \prod_{i=1}^{N} \mathrm{d H}_{ii}\prod_{i<j}\mathrm{d Re H}_{ij}\; \mathrm{d Im H}_{ij}$ is Lebesgue measure and the normalization constant $Z^{\JUE}$ is computed as
\begin{equation*}
  \label{eq:Z-JUE}
  Z^{\JUE}=\int (\det H)^{\alpha} (\det(I-H))^{\beta} \mathrm{d H}= \pi^{\frac{n(n-1)}{2}}\prod_{k}\frac{\Gamma(\alpha+k)\Gamma(\beta+k)}{\Gamma(\alpha+\beta+2n+1-k)}.
\end{equation*}
If $\alpha$ and $\beta$ are integers, so that $M_{\alpha}=\alpha+N$ and $M_{\beta}=\beta+N$ are integers, this probability measure describes the distribution of the matrix that is constructed from  Wishart random matrices $W_{A}=A^{*}A$, $W_{B}=B^{*}B$ with $A$ being $M_{\alpha}\times N$ and $B$ being $M_{\beta}\times N$ complex random matrices with i.i.d. Gaussian entries as
\begin{equation*}
  \label{eq:JUE-from-Wishart}
  H=(W_{A}+W_{B})^{-\frac{1}{2}} W_{A}(W_{A}+W_{B})^{-\frac{1}{2}}.
\end{equation*}
Integrating over $U(N)$, we obtain the probability measure for the distribution of the eigenvalues $\lambda_{1},\dots,\lambda_{N}$
\begin{equation}
  \label{eq:JUE-eigenvalues-distribution}
  \dmu^{\JUE}_{\alpha,\beta}(\lambda_{1},\dots,\lambda_{N})=\frac{1}{Z^{\JUE}}\prod_{i<j} (\lambda_{i}-\lambda_{j})^{2}\prod_{i}\lambda_{i}^{\alpha}(1-\lambda_{i})^{\beta}\mathrm{d}\lambda_{i}.
\end{equation}
Introduce the empirical measure with the density $\sigma_{N}(x)=\frac{1}{N}\sum_{i=1}^{N}\delta(x-\lambda_{i})$. Then in the large $N$ limit such that $N,M_{\alpha},M_{\beta}$ go to infinity with the same rate, so that  $\lim\frac{M_{\alpha}}{N}= c_{\alpha}$ and $\lim\frac{M_{\beta}}{N}=c_{\beta}$ with $c_{\alpha}+c_{\beta}>1$, the empirical measure converges weakly in probability  $\sigma_{N}(x)\dx\to\sigma(x)\dx$ to measure with the density
\begin{equation}
  \label{eq:JUE-limiting-law}
  \begin{aligned}
      &\sigma(x)=\frac{\sqrt{(\lambda_{+}-x)(x-\lambda_{-})}}{2\pi x(1-x)}\mathds{1}_{(\lambda_{-},\lambda_{+})}(x)+ \\
      & +\max(0,1-c_{\alpha})\delta(x)+\max(0,1-c_{\beta})\delta(x-1),
  \end{aligned}
\end{equation}
where the boundaries of the continuous part of the density are
\begin{equation*}
  \label{eq:JUE-continuous-domain}
  \lambda_{\pm}=\frac{c_{\alpha}(c_{\alpha}+c_{\beta}-1) + c_\beta\pm 2\sqrt{c_{\alpha}c_{\beta}(c_{\alpha}+c_{\beta}-1)}}{(c_{\alpha}+c_{\beta})^{2}}. 
\end{equation*}
This density generalizes Wigner's semicircle law and Marchenko--Pastur law. This result was obtained in \cite{capitaine2004asymptotic} and, for the sum of Wishart matrices in \cite{pastur2000law},  see also \cite{ledoux2004differential,collins2005product} and Proposition 3.6.3 in \cite{forrester2010log}. A plot of this density and a random samples for different $N$ is presented in Figure~\ref{fig:JUE-limit-density}. To recover this result we, for example, can write \eqref{eq:JUE-eigenvalues-distribution} in the exponential form as
\begin{equation*}
  \label{eq:JUE-eigenvalues-exponential}
  \dmu^{\JUE}_{\alpha,\beta}(\sigma_{N})=\exp(-N^{2}J[\sigma_{N}])
\end{equation*}
in terms of $\sigma_{N}$ and minimize the functional
\begin{equation*}
  \label{eq:JUE-functional}
  J[\sigma]=\int_{0}^{1}\int_{0}^{1}\sigma(x)\sigma(y)\ln|x-y|^{-1}\dx \dy+\int_{0}^{1}\sigma(x) V(x)\dx
\end{equation*}
with $V(x)=(1-c_{\alpha})\ln x+(1-c_{\beta})\ln (1-x)$ by solving the variational problem for the density $\sigma$. Taking variation of $J$ with respect to $\sigma$ and differentiating with respect to $x$ we obtain the following equation:
\begin{equation*}
  \label{eq:electrostatic-equilibrium}
  -\mathrm{v.p.}\int_{\mathrm{supp}\sigma}\frac{\sigma(y)\dy}{y-x}+\frac{1}{2}V'(x)=0.
\end{equation*}
It can be interpreted as the electrostatic equilibrium condition for the charged two-dimensional particles on real line. The density $\sigma$ is then a density of charged particles held in equilibrium by the external field $V'$. Denote by $G_{\sigma}(z)$ the electric field strength $G_{\sigma}(z)=-\int\frac{\sigma(y)\dy}{y-z}$, then it is analytic on $\CC\setminus \supp\sigma$ and $\lim_{\varepsilon\to 0}G_{\sigma}(x\pm i\varepsilon)=\mp i \pi \sigma(x)+ \frac{1}{2}V'(x)$. Therefore $G_{\sigma}(z)$ satisfies Riemann--Hilbert problem. Assuming the support of $\sigma$ to be single interval $[a,b]$ we can write the solution in the form  \cite[equation (1.67)]{forrester2010log} or \cite[Chapter 6]{deift2000orthogonal}
\begin{equation*}
  \label{eq:RH-solution}
  \sigma(x)=\frac{1}{\pi^{2}}\sqrt{(x-a)(x-b)}\int_{a}^{b}\frac{V'(x)-V'(y)}{x-y}\frac{\dy}{\sqrt{(b-y)(y-a)}},
\end{equation*}
with $a, b$ obtained from equations
\begin{equation*}
  \label{eq:support}
  \int_{a}^{b}\frac{V'(y)\dy}{\sqrt{(b-y)(y-a)}}=0,\quad \frac{1}{\pi}\int_{a}^{b}\frac{y V'(y)\dy}{\sqrt{(b-y)(y-a)}}=1.
\end{equation*}
This solution, after substitution of the potential $V(x)$, can be simplified to the continuous part of \eqref{eq:JUE-limiting-law}. 
\begin{figure}[h]
  \centering
  \includegraphics[width=0.9\linewidth]{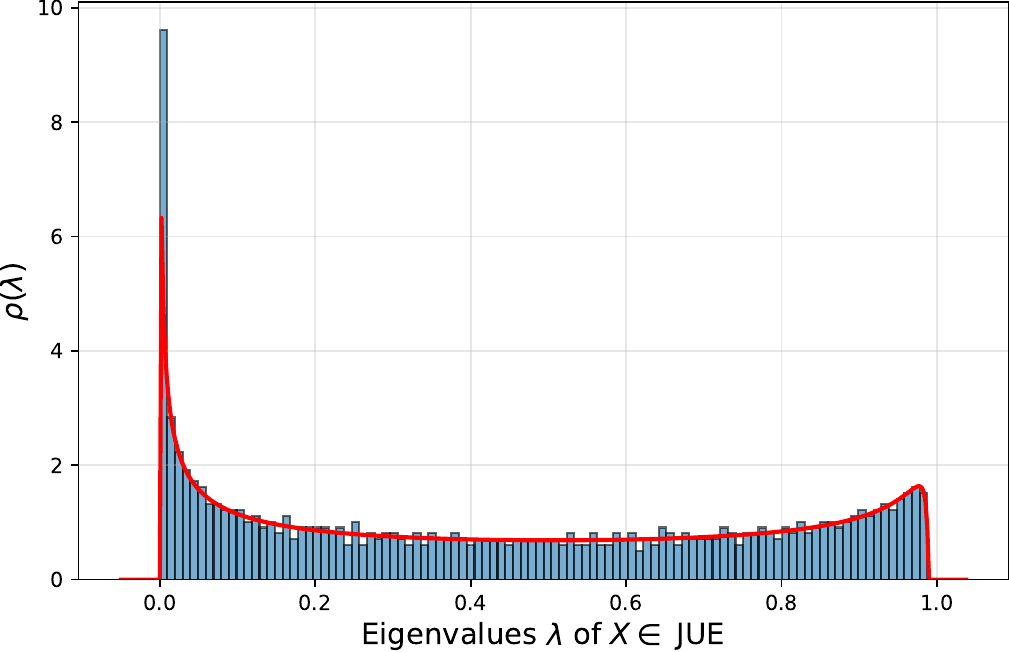}
  \caption{Histogram of JUE eigenvalues for $N=1000$ and $c_{\alpha}=0.95$, $c_{\beta}=1.235$ and limiting law (red). Note the discrete load at zero due to $c_\alpha < 1$. }
  \label{fig:JUE-limit-density}
\end{figure}

The limit shape for the random Young diagrams with respect to the measure~\eqref{eq:skew-howe-probability-measure} can be obtained in a similar way, see for example \cite{nazarov2024skew}. First, use Weyl dimension formula to write the measure as a product of Vandermonde determinant and a binomial coefficient
\begin{equation*}
  \label{eq:skew-howe-measure-formula}
  \mu^{\Howe}(\lambda)=\frac{1}{Z^{\Howe}}\prod_{i<j}(\lambda_{i}-\lambda_{j}+j-i)^{2} \prod_{i=1}^{n}\binom{n+k-1}{\lambda_{i}+n-i}.
\end{equation*}
We omit the expression for the normalization constant $Z^{\Howe}$ as it is not needed here. 
Then introduce the empirical particle density for the diagram $\lambda$ as $\rho_{n}(x)=\frac{1}{n}\sum_{i}\delta(x-\frac{\lambda_{i}-i}{n})$, 
rewrite the measure as the exponential of a functional $J[\rho_{n}]$ and apply Stirling approximation for $n,k,p\to\infty$ with $\lim\frac{k}{n}=c$, $\frac{p}{nk-p}=\alpha^{\Howe}$ to get the potential
\begin{equation*}
  \label{eq:Howe-potential}
  V(x)=(c-x)\ln(c-x)+(x+1)\ln(x+1)+\zeta x,
\end{equation*}
with the term $\zeta x$ coming from the condition $\lambda\vdash p$ and $\zeta$ being the Lagrange multiplier.

Minimizing the functional the same way as described above, we get Riemann--Hilbert problem for $G_{\rho}(z)=-\int\frac{\rho(x)\dx}{x-z}$. Solving it with the additional condition we obtain $\zeta=\ln \alpha^{\Howe}$ and 
\begin{equation*}
  \label{eq:Howe-limit-density}
  \rho(x)=\frac{1}{\pi}\arccos\left(\frac{\alpha^{\Howe}(c-1)+x(1-\alpha^{\Howe})}{2\sqrt{\alpha^{\Howe} (c-x)(x+1)}}\right),
\end{equation*}
for $x$ in the support of $\rho$, where $\operatorname{supp} \rho(x)=\left[t_{-},t_{+}\right]\subset [-1,c]$ with
  \begin{equation*}
    \label{eq:support-interval}
    t_{\pm}=\frac{\alpha^{\Howe}(c-1)\pm2\sqrt{\alpha^{\Howe}c}}{\alpha^{\Howe}+1}.
  \end{equation*}
  See Figure~\ref{fig:JUE-limit-density} for a rescaled plot of limit shape $\Omega(x)=1+\int_{-1}^{x}(1-2\rho(t))\mathrm{dt}$ obtained by integrating the density $\rho(x)$.
  
  Markov--Krein correspondence \eqref{eq:markov-krein} can be seen as the relation $G_{\sigma}(u)=\frac{1}{u+1}\exp G_{\rho}(u)$ between $G_{\sigma}$ and $G_{\rho}$. Computing the transition measure $\sigma$ from this relation we obtain
\begin{multline*}
  \label{eq:transition-measure}
  \sigma^{\Howe}(x)=\max\left(0,\frac{\alpha^{\Howe}c-1}{\alpha^{\Howe}(c+1)}\right) \delta(x+1)+\max\left(0,\frac{\alpha^{\Howe}-c}{\alpha^{\Howe}(c+1)}\right) \delta(c-x)+\\+\frac{1+\alpha^{\Howe}}{2\pi \alpha^{\Howe}}\frac{\sqrt{(x-t_{-})(t_{+}-x)}}{(c-x)(x+1)}.
\end{multline*}
By choosing $c=\frac{c_{\beta}}{c_{\alpha}}$ and $\alpha^{\Howe}=\frac{1}{c_{\alpha}+c_{\beta}-1}$ and rescaling the interval $(-1,c)$ to $(0,1)$ we recover \eqref{eq:JUE-limiting-law}.

\section{Random diagrams and counting of Riemann surfaces}
\label{sec:rand-diagr-count}

Consider now the finite $n,k,p$. The space $\bigwedge^{p}(\CC^{n}\otimes\CC^{k})$ decomposes into $\binom{nk}{p}$ sign representations of the symmetric group $S_{p}$ and so an introduction of transition probability by formula \eqref{eq:transition-probability} does not seem to be related to this case. But we can consider the exterior power as a module over the center of the symmetric group algebra $Z\CC[S_{p}]$, as was demonstrated in \cite{panova2018skew}. In this section we will compute the characters of powers of Young--Jucys--Murphy elements by using Stanley's character formula \cite{stanley2004irreducible}.

The group algebra $\CC[S_{p}]$ consists of linear combinations $\sum_{w\in S_{p}}a_{w}w$ with complex coefficients $a_{w}\in\CC$. Center of the group algebra $Z\CC[S_{p}]$ is spanned by the basis elements $\mathcal{C}_{\lambda}=\sum_{w\in[\lambda]}w$ for $\lambda\vdash p$, where $[\lambda]$ is the conjugacy class of permutations that have cycle type $\lambda$. Young--Jucys--Murphy elements are defined as sums of transpositions $\mathcal{J}_{1}=0$, $\mathcal{J}_{a}=(1a)+(2a)+\dots+(a-1 \ a)$, $a=2,\dots,p$. These elements are not central, but commute among themselves. Symmetric polynomials of $p$ variables evaluated at $\mathcal{J}_{1},\dots, \mathcal{J}_{p}$ generate the center $Z\CC[S_{p}]$. Moreover $\prod_{a=1}^{p}(1+\varepsilon \mathcal{J}_{a})=\sum_{\lambda\vdash p}\varepsilon^{p-\ell(\lambda)}\mathcal{C}_{\lambda}$, where $\ell(\lambda)$ is the number of rows in $\lambda$. 

Sniady and Panova in \cite{panova2018skew} defined the character of $Z\CC[S_{p}]$-module $V_{p}:=\bigwedge^{p}(\CC^{n}\otimes \CC^{k})$ as a function $\chi_{V_{p}}:S_{p}\to\CC$:
\begin{equation*}
  \label{eq:character-definition}
  \chi_{V_{p}}(w):=\frac{1}{\dim V_{p}}\tr [\Pi \rho_{w}].
\end{equation*}
Here $\rho_{w}$ is the action of $w$ on $W_{p}=(\CC^{n}\otimes\CC^{k})^{\otimes p}$, which permutes the factors in the tensor product. One can obtain the space $V_p$ as a subspace of $W_p$ by applying the projection $\Pi_{p}=\frac{1}{p!}\sum_{w\in S_{p}}(-1)^{w}w\in Z\CC[S_{p}]$. This means that the character of any permutation from $S_p$ of cycle type $\lambda$ can be computed as the character of the corresponding basis element $C_\lambda$ of the center $Z\CC[S_{p}]$ divided by the cardinality of $[\lambda]$. Note that we can extend the character $\chi_{V_{p}}$ to $\CC[S_{p}]$ by linearity, and that this character is normalized so that $\chi_{V_p}(\text{Id}) = 1$. 

According to Lemma 3.2 in \cite{panova2018skew}, for any $p'>p$ the character $\chi_{V_{p}}$ is equal to the restriction of $\chi_{V_{p'}}$. Therefore, it is sufficient to consider only the character $\chi_{V_{nk}}$. Corollary 3.3 in \cite{panova2018skew} provides that the character $\chi_{V_{p}}$, as the restriction of the character $\chi_{V_{nk}}$, is equal to the restriction of the  $S_{nk}$ character  $\chi_{k^{n}}: S_{nk} \rightarrow \CC$ of the representation corresponding to the rectangular $n\times k$ Young diagram. To compute the latter we use the results of \cite{stanley2004irreducible}. Let $\mathcal{X}$ denote the unnormalized character of the symmetric group, satisfying $\mathcal{X}_\lambda(\mathrm{Id}) = \dim(\lambda)$. For a diagram $\lambda\vdash p$ consider the diagram $(\lambda,1^{nk-p}) \vdash nk$ obtained by padding $\lambda$ by $nk-p$ rows of unit length. For $\mu\vdash nk$, define 
\begin{equation*}
    \label{eq:normalized-character-Stanley}
    \hat{\chi}_{\mu}(\lambda,1^{nk-p})=\frac{(nk)!\mathcal{X}_{\mu}(\lambda,1^{nk-p})}{(nk-p)! \mathcal{X}_{\mu}(1^{nk})}.
\end{equation*}
For $\mu=k^{n}$, formula (7) in \cite{stanley2004irreducible} gives the following expression:
\begin{equation}
  \label{eq:normalized-character}
  \hat{\chi}_{k^{n}}(\lambda,1^{nk-p})=(-1)^{p}\sum_{\nu\vdash p}H_{\nu} s_{\nu}(1^{n})s_{\nu}(1^{-k})\mathcal{X}_{\nu}(\lambda),
\end{equation}
where $H_{\nu}=\prod_{\square_{ij}\in\nu}h(\square_{ij})$ is the product of hooks of diagram $\nu$, $s_{\nu}$ are Schur polynomials and $s_{\nu}(1^{n})=\dim V_{\GL_{n}}(\nu)$, $s_{\nu}(1^{-k})=(-1)^{|\nu|}s_{\nu'}(1^{k})=(-1)^{p}\dim V_{\GL_{k}}(\nu')$. Hook product can be written as $H_{\nu}=\frac{p!}{\dim \nu}$. From \eqref{eq:normalized-character} we then obtain:
\begin{equation*}
  \label{eq:character-computation} \frac{\mathcal{X}_{k^{n}}(\lambda,1^{nk-|\lambda|})}{\mathcal{X}_{k^{n}}(1^{nk})}=\frac{(nk-p)!p!}{(nk)!}\sum_{\nu\vdash p}\frac{1}{\dim \nu}\dim V_{\GL_{n}}(\nu)\dim V_{\GL_{k}}(\nu') \mathcal{X}_{\nu}(\lambda).
\end{equation*}
The left-hand side is exactly the restriction of the $S_{nk}$-character $\chi_{k^n}$ to $S_p$. Now to compute $\chi_{V_{p}}(\mathcal{J}_{p}^{m})$ we should expand $\mathcal{J}_{p}^{m}$ as a linear combination over diagrams $\lambda$ of at most $p$ boxes, then for each $\lambda$ compute the character. The final step is to take into account that eigenvalues of Young--Jucys--Murphy elements in irreducible $S_{p}$-representation corresponding to diagram $\nu$ are contents of the boxes that can be removed from $\nu$, so $\mathcal{X}_{\nu}(\mathcal{J}_{p}^{m})=\sum_{i}\dim (\nu-\square_{i}) (\nu_{i}-i)^{m}$, where the sum is over the rows where the last box $\square_{i}$ can be removed to produce Young diagram with one less box \cite{okounkov2000random}. The ratio $\delta_{i}(\nu)=\frac{\dim(\nu-\square_{i})}{\dim \nu}$ is a dual object to the transition probability \eqref{eq:transition-probability} and is called cotransition probability or decay rate. Therefore we obtain
\begin{equation*}
  \label{eq:Young-Jucys-Murphy-character}
  \chi_{V_{p}}(\mathcal{J}_{p}^{m})=\sum_{\nu\vdash p} \mu^{\Howe}(\nu) \sum_{i}\delta_{i}(\nu) (\nu_{i}-i)^{m},
\end{equation*}
where we have substituted the dimension of $V_{p}=\bigwedge^{p}(\CC^{n}\otimes \CC^{k})$, $\dim V_{p}=\binom{nk}{p}$.

This computation can be iterated for other Young--Jucys--Murphy elements $\mathcal{J}_{a}$, $a<p$, since the $\mathcal{J}_{p}$-eigenspace, that corresponds to the eigenvalue $(\nu_{i}-i)$, is the irreducible representation of the subgroup $S_{p-1}$ in $S_{p}$ that fixes $p\in\{1,\dots,p\}$ and corresponds to the diagram $\nu-\square_{i}$. Therefore the eigensubspaces of $\mathcal{J}_{p-1}$ in this eigenspace correspond to the corners of the diagram $\nu-\square_{i}$ and are irreducible over the subgroup $S_{p-2}$ that fixes $p$ and $p-1$. See \cite[equation (3.14)]{okounkov2000random}. 

On the other hand, any permutation $w_{\lambda}$ in $S_p$ of the cycle type $\lambda$ acts on every basis vector in $V_p$ by multiplication by $\operatorname{sign}(\lambda) = (-1)^{|\lambda| - \ell(\lambda)}$, so the corresponding character is:
\begin{equation*}
\label{eq:character-sign}
    \chi_{V_p}(w_\lambda) = (-1)^{|\lambda|-\ell(\lambda)}. 
\end{equation*}
The product $\prod_{i=1}^s \mathcal{J}_{p+1-i}^{\kappa_i}$ can be written as a sum of terms of the following form:
\begin{equation}
    \label{eq:transpositions-product}
    w(\tau) := (\tau_{11} p)\dots(\tau_{1\kappa_1}p)(\tau_{21} \ p-1)\dots(\tau_{2\kappa_2} \ p-1)\dots (\tau_{s\kappa_s} \ p+1-s).
\end{equation}
Every transposition changes the sign of the final permutation $w(\tau)$. This means that every term contributes $(-1)^{\kappa_1 + \kappa_2 + \dots + \kappa_s} = (-1)^{|\kappa|}$ to the value of $\chi_{V_p}$ applied to $\prod_{i=1}^s \mathcal{J}_{p+1-i}^{\kappa_i}$:
\begin{equation}
    \label{eq:character-sign-YJM-cardinality}
    \chi_{V_p}\left(\prod_{i=1}^s \mathcal{J}_{p+1-i}^{\kappa_i}\right) = (-1)^{|\kappa|}|\{\tau\}|,
\end{equation}
where $\{\tau\}$ is the set of all possible choices of individual $\tau_{ij}$. Each $\tau_{ij}$ takes values from $1$ to $p-i$, so the cardinality in \eqref{eq:character-sign-YJM-cardinality} can be easily calculated:
\begin{equation}
    \label{eq:character-YJM-exact}
     \chi_{V_p}\left(\prod_{i=1}^s \mathcal{J}_{p+1-i}^{\kappa_i}\right) = (-1)^{|\kappa|}(p-1)^{\kappa_1} (p-2)^{\kappa_2}\dots (p-s)^{\kappa_s}.
\end{equation}

The third expression for the character of Young--Jucys--Murphy elements can be obtained by using Stanley's formula \cite[Theorem 1] {stanley2004irreducible}. According to this formula, the character of the permutation $w(\tau) \in S_p$ is given by
\begin{equation*}
    \label{eq:Stanley's-formula}
    \frac{\mathcal{X}_{k^n}(w(\tau),1^{nk-p})}{\mathcal{X}(1^{nk})} = (-1)^p\frac{(nk-p)!}{(nk)!} \sum_{ w_\nu w_\mu w(\tau)= \text{Id}}n^{\ell(\mu)}(-k)^{\ell(\nu)},
\end{equation*}
where the sum runs over all $p!$ pairs $w_\mu,w_\nu \in S_p$ of arbitrary cycle types $\mu, \nu$ satisfying $w_\nu w_\mu w(\tau)= \text{Id}$. Summing over all $\tau$ and rearranging the final sum, we obtain 
\begin{equation*}
    \label{eq:character-YJM-topological}
    \begin{aligned}
        & \chi_{V_p}\left(\prod_{i=1}^s \mathcal{J}_{p+1-i}^{\kappa_i}\right) = (-1)^p\frac{(nk-p)!}{(nk)!} \sum_\tau \sum_{ w_\nu w_\mu w(\tau)= \text{Id}}n^{\ell(\mu)}(-k)^{\ell(\nu)} =\\
        & =  (-1)^p\frac{(nk-p)!}{(nk)!} \sum_{g \ge 0} n^{2-2g+|\kappa|}\sum_{\mu, \nu \vdash p} \left(-\frac{k}{n}\right)^{\ell(\nu)} \mathcal{H}_g(\kappa,\mu,\nu).
    \end{aligned}
\end{equation*}
Here $\mathcal{H}_g(\kappa, \mu, \nu)$ is defined to be the number of tuples $(w_\nu, w_\mu, w(\tau))$ of permutations in $S_p$ such that
\begin{enumerate}
    \item $2-2g = \ell(\mu) + \ell(\nu) - |\kappa|$,
    \item $w_\nu \in [\nu], w_\mu \in [\mu]$,
    \item $w(\tau)$ is defined in \eqref{eq:transpositions-product},
    \item $w_\nu w_\mu w(\tau) = \text{Id}$.
\end{enumerate}

The relation $w_\nu w_\mu w(\tau) = \text{Id}$, where $w(\tau)$ is defined in \eqref{eq:transpositions-product}, can be seen as a monodromy representation for the Riemann covering of a sphere with $|\kappa|$ branch points having simple ramification and two additional points with ramification profiles $\mu, \nu$. We can assume all branching points to be distinct from $0$ and $\infty$ and all loops to start at $0$ and encircle the branching points. The product of all loops around branching points is contractible, and hence the product of corresponding monodromies must be equal to one.

The correlation functions for JUE are also expressed using Young--Jucys--Murphy elements. The generating functions for JUE correlators are hypergeometric Kadomtsev--Petviashvili tau functions \cite{harnad2015hypergeometric}. Below we adopt the notation from \cite{gisonni2021jacobi}. Consider the correlation functions of Jacobi Unitary Ensemble  $\left<\prod_{i=1}^{\ell(\kappa)}\tr H^{\kappa_{i}} \right>$, where we can assume $\kappa_{1}\geq \kappa_{2}\geq\dots\geq \kappa_{\ell(\kappa)}$ and $\kappa$ to be a Young diagram.  Denote by $Z_{N}(u_{1},u_{2},\dots)$ the formal generating function in infinitely many variables $(u_{1},u_{2},\dots)$ for the correlation functions
\begin{equation*}
  \label{eq:JUE-generating-function}
  \begin{aligned}
      & Z_{N}(u_{1},\dots)=\sum_{\lambda}\frac{1}{z_{\kappa}}\left< \prod_{j=1}^{\ell(\kappa)}\tr H^{\kappa_{j}}\right> u_{\kappa_{1}}\dots u_{\kappa_{j}}= \\
      & = \int\exp\left(\sum_{m\geq 1}\frac{u_{m}}{m}\tr H^{m}\right) \mathrm{dm}^{\JUE}_{N}(H),
  \end{aligned}
\end{equation*}
where $z_{\kappa}=\prod_{i\geq 1} i^{m_{i}}(m_{i}!)$ and $m_{i}$ is the number of rows in $\kappa$ of length $i$.

Consider the space of symmetric functions of infinitely many variables $x_{1},x_{2},\dots$. It can be seen as spanned by the power sums $p_{m}=\sum_{i}x_{i}^{m}$, by elementary symmetric polynomials, homogeneous symmetric polynomials or Schur polynomials $s_{\lambda}(x_{1},x_{2},\dots)$. 
Let  $\gamma(z)=\frac{(1+z)(1+\frac{z}{c_{\alpha}})}{1+\frac{z}{c_{\alpha}+c_{\beta}}}$ and $r_{\kappa}(\gamma,\varepsilon)=\prod_{\square_{ij}\in\kappa}\gamma(\varepsilon (j-i))$ . Then the hypergeometric tau function can be defined by
\begin{equation*}
  \label{eq:tau-function-schur-expansion}
  \tau_{\gamma}(\varepsilon;p_{1},p_{2},\dots)=\sum_{\kappa}\frac{\dim \kappa}{|\kappa|!}r_{\kappa}(\gamma,\varepsilon) s_{\kappa}(x_{1},x_{2},\dots). 
\end{equation*}
For any $p>0$ a multiparametric Hurwitz number $H_{\gamma}^{d}(\kappa)$ is defined as the coefficient of $\varepsilon^{d}\mathcal{C}_{\kappa}$ in class basis expansion of $\prod_{a=1}^{p}\gamma(\varepsilon\mathcal{J}_{a})$
\begin{equation*}
  \label{eq:multiparametric-Hurwitz-numbers}
  H_{\gamma}^{d}(\kappa)=\frac{1}{z_{\kappa}}[\varepsilon^{d} \mathcal{C}_{\kappa}]\prod_{a=1}^{p}\gamma(\varepsilon\mathcal{J}_{a}).
\end{equation*}
The hypergeometric tau function admits the expansion in terms of these Hurwitz numbers as 
\begin{equation*}
  \label{eq:tau-function-expansion}
  \tau_{\gamma}(\varepsilon;p_{1},p_{2}\dots)=\sum_{d\geq 1} \varepsilon^{d} \sum_{\kappa} H_{\gamma}^{d}(\kappa)\prod_{i=1}^{\ell(\kappa)}p_{\kappa_{i}}.
\end{equation*}
The generating function for the correlators is the hypergeometric tau function for $\varepsilon=\frac{1}{N}$ and $p_{m}=\left(\frac{c_{\alpha}}{c_{\alpha}+c_{\beta}}\right)^{m}N^{m}u_{m}$
\begin{equation*}
  \label{eq:generating-function-is-tau}
  Z_{N}(u_{1},u_{2},\dots)=\tau_{\gamma}\left(\frac{1}{N};\left(\frac{c_{\alpha}}{c_{\alpha}+c_{\beta}}\right)N u_{1},\dots,\left(\frac{c_{\alpha}}{c_{\alpha}+c_{\beta}}\right)^{m}N^{m}u_{m},\dots\right).
\end{equation*}
Using this result it is possible to write genus expansion for the correlation functions $\left< \prod_{j=1}^{\ell(\kappa)}\tr H^{\kappa_{j}}\right>$ as
\begin{equation}
  \label{eq:genus-expansion}
  N^{\ell(\kappa)}\frac{|\kappa|!}{z_{\kappa}}\left<\prod_{j=1}^{\ell(\kappa)}\tr H^{\kappa_{j}}\right>=\sum_{g\geq 0}\frac{(-1)^{|\kappa|}}{N^{2g-2}}\sum_{\mu,\nu\vdash|\kappa|}\frac{c_{\alpha}^{\ell(\nu)} h_{g}(\kappa,\mu,\nu)}{(-c_{\alpha}-c_{\beta})^{\ell(\mu)+\ell(\nu)+\ell(\kappa)+2g-2}},
\end{equation}
where $h_{g}(\kappa,\mu,\nu)$ is called triple monotone Hurwitz number. It is the number of tuples $(w_\mu,w_\nu,\tau_{1},\dots,\tau_{r})$, with
\begin{enumerate}
\item $r=2g-2-|\kappa|+\ell(\kappa)+\ell(\mu)+\ell(\nu)$,
\item $w_\mu\in [\mu]$, $w_\nu\in [\nu]$,
\item transpositions $\tau_{i}=(a_{i},b_{i})$, $a_{i}<b_{i}$ are ordered $b_{1}\leq \dots\leq b_{r}$,
\item $w_\mu w_\nu\tau_{1}\dots\tau_{r}\in[\kappa]$.
\end{enumerate}
Triple monotone Hurwitz numbers count Riemann surfaces of degree $d$ with three branching points with branching profiles $\kappa,\mu,\nu$ and simple branching at $r$ points. They are connected to multiparametric Hurwitz numbers by the formula
\begin{equation*}
  \label{eq:triple-monotone-and-multiparametric-Hurwitz}
  H^{2g-2+|\kappa|+\ell(\kappa)}_{\gamma}(\kappa)=\frac{1}{|\kappa|!}\sum_{\mu,\nu\vdash|\kappa|}c_{\alpha}^{\ell(\nu)-|\kappa|}(-c_{\alpha}-c_{\beta})^{|\kappa|+2-2g-\ell(\mu)-\ell(\nu)-\ell(\kappa)}h_{g}(\kappa,\mu,\nu).
\end{equation*}
Correlators of large powers of $H$ in JUE are dominated by the eigenvalues near the edge of the spectrum $\lambda_{+}$.  In the asymptotic regime $\kappa_{i}\sim \xi_{i}N^{\frac{2}{3}}$ (soft edge) correlators are described by the Airy kernel \cite{forrester2010log}. The fluctuations of the largest eigenvalue obey Tracy--Widom distribution \cite{TW94}. In this regime the terms of genus expansion \eqref{eq:genus-expansion} also depend on $N$ through $\kappa$, therefore we need more than the asymptotic of a single Hurwitz number $h_{g}(\kappa,\mu,\nu)$. The JUE version of Hurwitz measure \cite{okounkov2006gromov,okounkov2009gromov} needs to be introduced and studied asymptotically.

Similarly the correlators of particle positions near the edge of the limit shape for $\mu^{\Howe}$ are described by the Airy kernel, and the fluctuations of the first particle position or the first row length obey Tracy--Widom distribution \cite[Section 4.2]{betea2025limit}. Here we need to consider the regime $\kappa_{i}\sim \xi_{i}n^{\frac{2}{3}}\sim \xi_{i}p^{\frac{1}{3}}$. In this regime we can neglect the permutations between sheets numbered $p,p-1,\dots,p-s+1$ in \eqref{eq:transpositions-product} similar to \cite[equation (3.15)]{okounkov2000random} as it does not change the leading asymptotic of order $p^{p^{1/3}}$ in \eqref{eq:character-YJM-exact}. For random diagram $\lambda$ and finite $i$ in this regime we have $\lambda_{i}-i\approx \lambda_{i}$.  Therefore we can conjecture that the joint distribution of random variables $n^{\frac{2}{3}}\left(\frac{\lambda_{i}}{t_{+}n}-1\right)$ corresponding to row lengths in random Young diagram with respect to $\mu^{\Howe}$ is identical to the join distribution of random variables $N^{\frac{2}{3}}\left(\frac{\lambda_{i}}{\lambda_{+}}-1\right)$ corresponding to the largest JUE eienvalues. We expect that leading asymptotic contributions to these distributions are given by the same set of Riemann surfaces the same way as for GUE/Plancherel correspondence in \cite{okounkov2000random}.

\section*{Conclusion and conjectures}
\label{sec:conjectures}
We have demonstrated that limit shape of random Young diagrams with respect to skew Howe measure $\mu^{\Howe}$ in the limit $n,k,p\to \infty$ with $p\sim nk$ is related to the limiting distribution of the JUE eigenvalues when $\alpha,\beta,N\to\infty$ with $\alpha\sim N$ and $\beta\sim N$.

We have also computed the characters of Young--Jucys--Murphy elements in $\bigwedge^{p}(\CC^{n}\otimes\CC^{k})$ and discussed their connection to Riemann surface counting. 

Let us state three conjectures that give the directions for the future work.
\begin{enumerate}
\item For random Young diagram $\lambda$ with respect to $\mu^{\Howe}$ in the asymptotic regime $n,k,p\to\infty$, $p\sim nk$ introduce the random variables $x_{i}=n^{\frac{2}{3}}\left(\frac{\lambda_{i}}{t_{+}n}-1\right)$ for $i=1,\dots, s$ for any finite $s$ and Laplace transform $\hat{x}(\xi)=\sum_{i=1}^{\infty}\exp(\xi x_{i})$.
  
  For random matrix $H$ from JUE with eigenvalues $\lambda_{1}\geq \lambda_{2}\geq\dots\geq \lambda_{N}$ in the asymptotic regime $\alpha,\beta,N\to\infty$, $\alpha\sim N$, $\beta\sim N$ introduce the random variables $y_{i}=N^{\frac{2}{3}}\left(\frac{\lambda_{i}}{\lambda_{+}}-1\right)$ for $i=1,\dots, s$ for any finite $s$ and Laplace transform $\hat{y}(\xi)=\sum_{i=1}^{\infty}\exp(\xi y_{i})$. 
  We conjecture that all mixed moments of the random variables $\hat{x}(\xi)$ and $\hat{y}(\xi)$ exist and are identical to each other and therefore the joint distribution of random variables $x_{1},\dots,x_{s}$ is identical to the joint distribution of $y_{1},\dots,y_{s}$ for any $s$.
\item The generating function for JUE correlators is tau function of an integrable hierarchy and admits a free-fermionic formulation \cite{okounkov2000todaequationshurwitznumbers, orlov2002matrixintegralssymmetricfunctions}, see also \cite{natanzon2020hurwitz}. In particular, there exists an integral representation for tau function \cite[Proposition 14]{orlov2002matrixintegralssymmetricfunctions}. There is a free-fermionic construction that gives an integral representation for the correlation kernel for $\mu^{\Howe}$, see for example \cite[Theorem 1.1]{betea2025limit}. We conjecture that these integral representations are asymptotically connected.
\begin{figure}[h]
  \centering
  \includegraphics[width=1.0\linewidth]{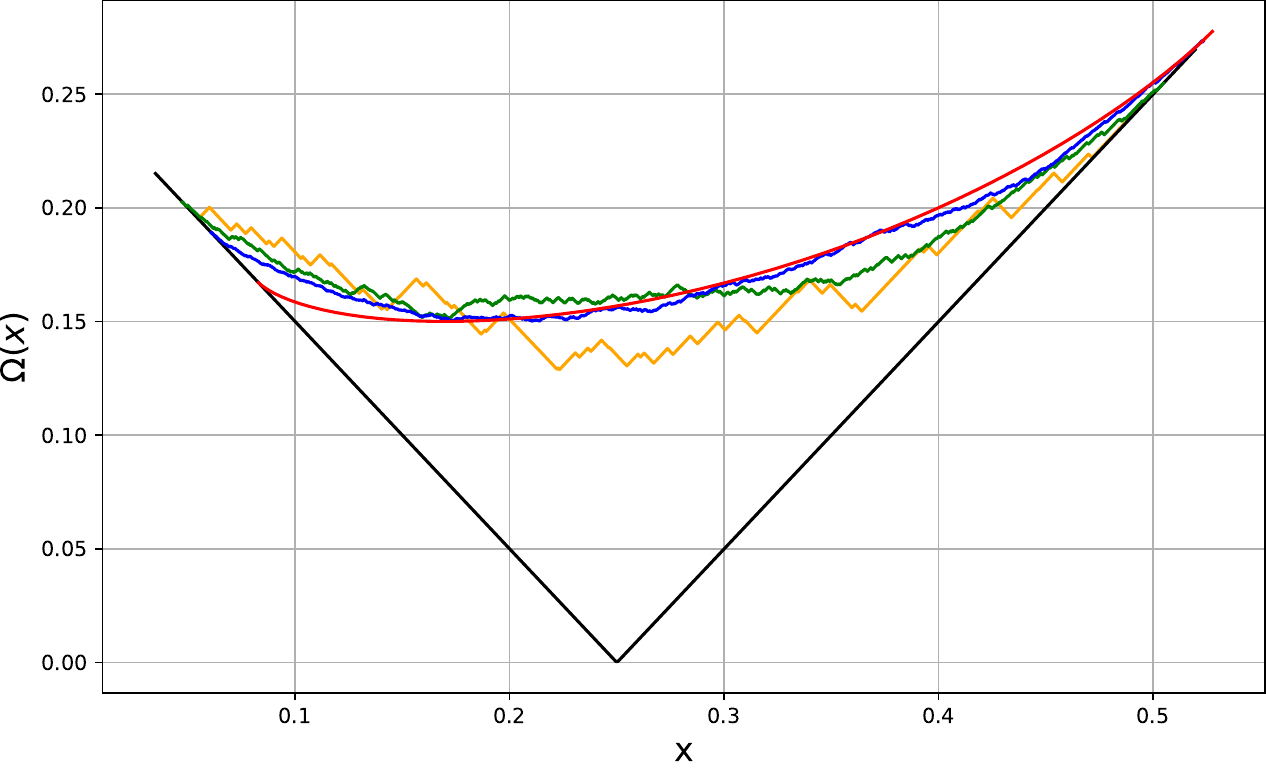}
  \caption{Piecewise-linear function constructed by interlacing eigenvalue of JUE random matrix $H$ and matrix $\hat{H}$, obtained from $H$ by removing last row and column, with $c_\alpha = 1.6, \ c_\beta = 4.8$ and $N = 50$ (\textcolor{YellowOrange}{orange}), $N = 500$ (\textcolor{OliveGreen}{green}), $N = 1500$ (\textcolor{blue}{blue}), and a rescaled limit shape of random Young diagrams for the measure $\mu^{\Howe}$ (\textcolor{red}{red}).}
  \label{fig:JUE-limit-density}
\end{figure}

\item   Construct a random piecewise-linear function $f$ by interlacing eigenvalues of JUE random matrix $H$ with the eigenvalues of random matrix $\hat{H}$, obtained by removing last row and last column of $H$ and then taking eigenvalues of $H$ to be positions of minima and eigenvalues of $\hat{H}$ to be positions of maxima and requiring $f'(x)=1$ between a minimum and the next maximum and $f'(x)=-1$ between a maximum and the next minimum. We conjecture that in the limit $\alpha,\beta,N\to\infty$ with $\alpha\sim N$, $\beta\sim N$ such random function converges to the limit shape $\Omega$ for $\mu^{\Howe}$ in supremum norm. In Figure~\ref{fig:JUE-limit-density} we present results of numerical simulations that support this conjecture. 
\end{enumerate}

\bibliographystyle{zapiski}
\makeatletter
\renewcommand\@biblabel[1]{#1.\hfill}
\makeatother
\bibliography{literature}{}

\end{document}